\documentclass[12pt]{article}

\usepackage{amsmath,amssymb}
\usepackage{latexsym}
\usepackage{graphicx}
\usepackage{bm}
\usepackage{hyperref}

\newcommand{\R}{\mathbb{R}}

\date{}

\begin{document}


\centerline{}

\centerline{}

\centerline {\Large{\bf  Global existence and uniqueness of the solution 
 }}

\centerline{\Large{\bf to a nonlinear
 parabolic equation}}
$$ $$
\centerline{\Large{\bf Alexander G. Ramm}}

\centerline{}

\centerline{Mathematics Department, Kansas State University,}

\centerline{Manhattan, KS 66506-2602, USA}

\renewcommand{\thefootnote}{\fnsymbol{footnote}}
\footnotetext{Corresponding author: Email: ramm@ksu.edu}

\newtheorem{Theorem}{\quad Theorem}[section]

\newtheorem{Definition}[Theorem]{\quad Definition}

\newtheorem{Corollary}[Theorem]{\quad Corollary}

\newtheorem{Lemma}[Theorem]{\quad Lemma}

\newtheorem{Example}[Theorem]{\quad Example}

\newtheorem{Remark}[Theorem]{\quad Remark}

\newtheorem{thm}{Theorem}[section]
\newtheorem{cor}[section]{Corollary}
\newtheorem{lem}[section]{Lemma}
\newtheorem{dfn}[section]{Definition}
\newtheorem{rem}[section]{Remark}

\newcommand{\bee}{\begin{equation*}}
\newcommand{\eee}{\end{equation*}}
\newcommand{\be}{\begin{equation}}
\newcommand{\ee}{\end{equation}}
\newcommand{\ba}{\begin{align}}
\newcommand{\ea}{\end{align}}
\newcommand{\pn}{\par\noindent}
\newcommand{\RRR}{\mathbb{R}^3}

\begin{abstract} \noindent Consider the equation
$$
u'(t)-\Delta u+|u|^\rho u=0, \quad u(0)=u_0(x), 
(1),
$$
where $ u':=\frac {du}{dt}$, $ \rho=const >0, $
$x\in \R^3$, $t>0$.

Assume that $u_0$ is a smooth and decaying function,
$$\|u_0\|\:=\sup_{x\in \R^3, t\in \R_+}
|u(x,t)|.$$
It is proved that problem (1) has a unique
global solution and this solution satisfies the following estimate
$$\|u(x,t)\|<c, $$
where $c>0$ does not depend on $x,t$.
\end{abstract}

35K55.

\section{Introduction}
Let 
\be\label{eq1}
u'-\Delta u +|u|^{\rho} u=0, \quad u(0)=u_0; \quad u':=\frac {du}{dt}, 
\ee
where $\rho>0$, $t\in \R_+=[0,\infty)$, $x\in \R^3$,
$X$ is a Banach space of real-valued functions with the norm
 $\|u(x,t) \|:=\sup_{x\in \R^3, t\in R_+}|u(x,t)|$.   
We assume that 
\be\label{eq2}
\|u\|\le c.
\ee
We say that $u$ is a global solution to \eqref{eq1}
if $u$ exists $\forall t\ge 0$.

Our result is formulated in Theorem 1. Our method 
is simple and differs from the published results,
see \cite{1}, \cite{2} and references there.
 
The novel points in this work  are: 

a) There is no restriction on the upper bound of 
$\rho$. 

In \cite{1}, (Section 1.1) a nonlinear hyperbolic equation with the same non-linearity is studied
in a bounded domain, uniqueness of the solution is proved only for $\rho\le 2/(n-2)$, and existence is 
proved by a different method.
 The contraction mapping
theorem is not used. 

In \cite{2} the quasi-linear problems for parabolic equations are studied in Chapter 5 in a bounded domain and under the assumptions different from ours.
There are many papers and books on non-linear
problems for parabolic equations (see the bibliography in \cite{1}, \cite{2}.
 
b) Existence of the global solution is proved.

c) Method of the proof differs from the methods in 
the cited literature. 

Our result is formulated in Theorem 1:

{\bf Theorem 1}. {\em Problem \eqref{eq1} has a unique global solution in $X$ for any $u_0\in X$.}

\section{Proofs}

Let $g(x,t)=\frac{e^{-|x|^2}}{(4\pi t)^{3/2}}$. If $u$
solves \eqref{eq1} then
\be\label{eq3}
 \begin{split}
u(t)=-\int_0^td\tau \int g(x-y, t-\tau)|u|^\rho u dy
+\\
\int g(x-y,t)u_0(y)dy:=A(u)+F:=Q(u),
 \end{split}
\ee
where $\int:=\int_{\R^3}$. Let $X$ be the Banach space of continuous in $\R^3 \times R_+$ functions, $\R_+:=[0, \infty)$, $\|u\|:=\max_{x\in \R^3, t\in [0,T]} |u(x,t)|$. If $\|u\|\le R$ then $\|A(u)\|\le T R^{\rho +1}$, where the identity $\int g(x-y, t-\tau)dy=1$ was used. From \eqref{eq3} one gets
\be\label{eq4}
\|u\|\le T\|u\|^{\rho +1} +\|F\|.
\ee
Thus, $Q$ maps the ball $B(R)=\{u: \|u\|\le R\}$
into itself if $T$ is such that 
\be\label{eq5} 
TR^{\rho +1} +\|F\|\le R.
\ee
The $Q$ is a contraction on $B(R)$ if
$$\|Q(u)-Q(v)\|\le T(\rho+1)R^\rho \|u-v\|\le q\|u-v\|, \quad 0<q<1.$$
Thus, if
\be\label{eq6}
T(\rho+1)R^\rho \le q<1,
\ee
then $Q$ is a contraction in $B(R)$ in the Banach space  $X_T$ with the norm $\|\cdot\|$, $t\in [0,T]$.
We use the same notations for the norms in $X_T$ and in $X_{\infty}$.

We have proved that  

{\em For $T$ satisfying \eqref{eq5}- \eqref{eq6} there exists and is unique the solution to \eqref{eq1}, and this solution can be
obtained from \eqref{eq3} by iterations.}

The problem now is: 

{\em Does this solution exist and is unique on $R_+$?} 

 From our proof it follows that if the solution exists
 and is unique in $X_T$, then the solution exists and is unique in $X_{T_1}$ for some $T_1>T$. 

To prove that the solution $u(x,t)$ to \eqref{eq1}
exists on $\R_+$, assume the contrary: this solution
does not exist on any interval $[0, T_1),$ $T_1>T$,
where $T$ is the maximal interval of the existence of the continuous solution. Then $\lim_{t\to T-0}u(x,t)=
\infty$, because otherwise there is a sequence $t_n\to T-0$ such that $u(x,t_n)\to u(x,T)$ and one may
construct the solution defined on $[T,T_1]$, $T_1>T$, by using the local existence and uniqueness of the solution to \eqref{eq1} with the initial value $u(x,T)$  for $t\in [T,T_1]$. This contradicts the assumption that $T$ is the maximal interval of the existence of the continuous solution $u$. 

Thus, if $T<\infty$ then one has $\lim_{t\to T-0}u(x,t)=\infty$. Let us prove that this also leads to a contradiction. Then we have to conclude that $T=\infty$ and Theorem 1 is proved.

We need some estimates. Multiply \eqref{eq1} by $u$, 
integrate over $\R^3$ with respect to $x$, and then
integrate by parts the second term. The result is: 
\be\label{eq7}
0.5\frac{d N(u)}{dt} + N(grad u) + \int |u|^{\rho+2}dy=0,
\ee
where $N(u):=\int u^2dy$.
Integrate \eqref{eq7} with respect to time over $[0,T]$ and get
\be\label{eq8}
0.5 N(u(T)) + \int_0^T \left(N(grad u)+\int |u|^{\rho+2}dy\right)d\tau =0.5 N(u(0)).
\ee
Therefore,
\be\label{eq9}
N(u(t))\le c,\,\, \forall t\in [0,T], \quad \int_0^T N(grad u)d\tau\le c,\quad
\int_0^T d\tau \int |u|^{\rho +2}dy \le c,
\ee
where $c=0.5 N(u_0)$.  

{\bf Lemma 1.} {\em From \eqref{eq9} and \eqref{eq3} it follows that
\be\label{eq10}
\|u(x,t)\|< \infty \quad \forall t\in [0,T].
\ee
}

If \eqref{eq10} is proved then $T$ is not the maximal interval of the existence of the solution to \eqref{eq1}. This contradiction proves Theorem 1.

{\em Proof of Lemma 1.} One uses the H\"older inequality twice and gets
\be\label{eq11}
\begin{split}
\int_0^T d\tau\int g(x-y, t-\tau)|u|^{\rho +1}dy\le\\
\left(\int_0^T d\tau\int |u|^{\rho +2}dy\right)^{(\rho+1)/(\rho +2)} \left(\int_0^T d\tau\int g^{\rho +2}dy\right)^{1/(\rho +2)}\le \\
\left(\int_0^Td\tau \int|u|^{\rho+2}dy\right)^{(\rho+1)/(\rho +2)} \left(\int_0^T d\tau \int g^{\rho +2}dy\right)^{1/(\rho +2)}.
\end{split}
\ee
By the last inequality \eqref{eq9} it follows that
$\int_0^Td\tau \int|u|^{\rho+2}dy<c$ $\forall T>0$,
where $c>0$ is a constant independent of $T$.
The last integral in \eqref{eq11} is also bounded independently of $T$. It can be calculated analytically.

Thus, inequalities \eqref{eq11}, \eqref{eq9} and equation \eqref{eq3} imply 
\eqref{eq10}. 

Lemma 1 is proved. \hfill$\Box$ 
  
Therefore Theorem 1 is proved. \hfill$\Box$

The ideas related to the ones used in this paper
were developed and used in [3]--[5].


\end{document}